# Linear Equation in Finite Dimensional Algebra

Aleks Kleyn

ABSTRACT. In the paper I considered methods for solving equations of the form
$$a_{s\cdot 0}\, x\, a_{s\cdot 1} = b$$
in the algebra which is finite dimensional over the field.

## Contents

## 1. Preface

My recent research in the field of calculus ([4]) and affine geometry over division ring ([5]) has led me to the need to solve linear equations of the form

$$(1.1) \qquad a_{s\cdot 0}\, x\, a_{s\cdot 1} = b$$

or system of such equations over division ring. The main problem that interests me is to find the inverse transformation, because this is important operation for the transformation of the tensor, as well as for lifting and lowering the index in Euclidean space.

In this paper I explored the possibility of solving the simplest equation. In modern mathematics and physics, scientists consider structures where product can be nonassociative. So I wonder what results will remain, if I consider linear algebra in finite dimensional algebra $A$ over field.

Coefficients of the equation (1.1) belongs to tensor product $A \otimes A$. This allows us to apply previously developed methods to solve equation (1.1). In the paper, I consider two methods to solve equation (1.1).

First, I consider the algebra $A$ as vector space over the field $F$. This allows me replace equation (1.1) by the system of linear equations which we can solve.

Aleks_Kleyn@MailAPS.org.
http://sites.google.com/site/AleksKleyn/.
http://arxiv.org/a/kleyn_a_1.
http://AleksKleyn.blogspot.com/.

Immediately it becomes evident that the equation (1.1) may have one solution, infinitely many solutions or none.

Though the standard representation of linear equation is not always the most compact form of notation, it simplifies the construction. Second method to solve equation (1.1) is based on the standard representation of linear expression. When the equation has one solution, this method allows us to find an inverse map (the problem which is important for me in tensor calculus), as well as better understand what it means to the singular linear map (this is important to construct homology of algebra; in particular, to understand how could look the Cauchy-Riemann equations in algebra and what functions hold this equation).

Expression (4.8) as record of the identity transformation looks unusually. I saw such time of expression first time in equation [3]-(6.4.3). However, it was the beginning of my research in this field, and I was not ready at that time to understand all the depth I have seen.

## 2. Conventions

**Convention 2.1.** *In any expression where we use index I assume that this index may have internal structure. For instance, considering the algebra $A$ we enumerate coordinates of $a \in A$ relative to basis $\overline{\overline{e}}$ by an index $i$. This means that $a$ is a vector. However, if $a$ is matrix, then we need two indexes, one enumerates rows, another enumerates columns. In the case, when index has structure, we begin the index from symbol $\cdot$ in the corresponding position. For instance, if I consider the matrix $a_j^i$ as an element of a vector space, then I can write the element of matrix as $a_{\cdot j}^{\cdot i}$.* □

**Convention 2.2.** *I assume sum over index $s$ in expression like*

$$a_{s \cdot 0} x a_{s \cdot 1}$$

□

**Convention 2.3.** *For given field $F$, unless otherwise stated, we consider finite dimensional $F$-algebra.* □

**Convention 2.4.** *Let $A$ be free finite dimensional algebra. Considering expansion of element of algebra $A$ relative basis $\overline{\overline{e}}$ we use the same root letter to denote this element and its coordinates. However we do not use vector notation in algebra. In expression $a^2$, it is not clear whether this is component of expansion of element $a$ relative basis, or this is operation $a^2 = aa$. To make text clearer we use separate color for index of element of algebra. For instance,*

$$a = a^i \overline{e}_i$$

□

**Convention 2.5.** *If free finite dimensional algebra has unit, then we identify the vector of basis $\overline{e}_0$ with unit of algebra.* □

Without a doubt, the reader may have questions, comments, objections. I will appreciate any response.

## 3. Product in Algebra

Let $F$ be field and $A$ be finite dimensional algebra over field $F$. Let $\overline{\overline{e}}$ be the basis of algebra $A$ over field $F$. We define the product of basis vectors according to rule

$$\overline{e}_i \overline{e}_j = C_{ij}^k \overline{e}_k \tag{3.1}$$

where $C_{ij}^k$ are structural constants of algebra $A$ over field $F$. Since I do not assume operation to be neither commutative, nor associative, I do not have any constraint for structural constants.

From equation (3.1), it follows that we can get the product of $a, b \in A$ according to rule

$$(ab)^k = C_{ij}^k a^i b^j \tag{3.2}$$

## 4. Linear Equation in Associative Algebra

We write linear equation in associative algebra $A$ in the following form

$$a_{s\cdot 0}\, x\, a_{s\cdot 1} = b \tag{4.1}$$

According to the theorem [3]-9.1.9, we can write the equation (4.1) in standard form

$$a^{ij} \overline{e}_i x \overline{e}_j = b \tag{4.2}$$

$$a^{ij} = a_{s\cdot 0}{}^i\, a_{s\cdot 1}{}^j \quad a_{s\cdot 0} = a_{s\cdot 0}{}^i e_i \quad a_{s\cdot 1} = a_{s\cdot 1}{}^i e_i$$

According to the theorem [3]-9.1.10 equation (4.2) is equivalent to equation

$$x^i a_i^j = b^j \tag{4.3}$$

$$a_i^j = a^{kr} C_{ki}^p C_{pr}^j \tag{4.4}$$

**Theorem 4.1.** *If determinant*

$$\det \|a_i^j\| \neq 0 \tag{4.5}$$

*then equation (4.1) has only one solution.*

*If determinant equal $0$, then $F$-linear dependence of vector $b$ from vectors $a_i^j \overline{e}_j$ is condition of existence of solution. In this case, equation (4.1) has infinitely many solutions. Otherwise equation does not have solution.*

*Proof.* The statement of the theorem is corollary of the theory of linear equations over field. □

**Theorem 4.2.** *Let the equation (4.2) satisfies to condition (4.5). If we consider the equation (4.2) as transformation of algebra $A$, then we can write the inverse transformation in form*

$$x = c^{pq} \overline{e}_p b \overline{e}_q \tag{4.6}$$

*where components $c^{pq}$ satisfy to equation*

$$\delta_0^r \delta_0^s = a^{ij} c^{pq} C_{ip}^r C_{qj}^s \tag{4.7}$$

*Proof.* According to theorem 4.1, equation (4.2) has only one solution. Since $x$ is linear function of $b$, then we consider standard representation (4.6) of this function. From equations (4.2), (4.6), it follows that

(4.8) $$b = a^{ij}c^{pq}\overline{e}_i\overline{e}_p b\overline{e}_q\overline{e}_j = a^{ij}c^{pq}C_{ip}^r C_{qj}^s \overline{e}_r b\overline{e}_s$$

Since $b$ is arbitrary, then the equation (4.8) is record of identity transformation. Equation (4.7) is corollary of the equation (4.8) and the theorem [3]-9.1.15. □

## 5. Linear Equation in Nonassociative Algebra

We write linear equation in nonassociative algebra $A$ in the following form

(5.1) $$a_{1 \cdot s \cdot 0}\,(x\,a_{1 \cdot s \cdot 1}) + (a_{2 \cdot t \cdot 0}\,x)\,a_{2 \cdot t \cdot 1} = b$$

According to the theorem [3]-9.1.9, we can write the equation (5.1) in standard form

(5.2) $$\begin{cases} a_1^{ij}\overline{e}_i(x\overline{e}_j) + a_2^{ij}(\overline{e}_i x)\overline{e}_j = b \\ a_1^{ij} = a_{1\cdot s\cdot 0}^{\,i}\, a_{1\cdot s\cdot 1}^{\,j} \quad a_{1\cdot s\cdot 0} = a_{1\cdot s\cdot 0}^{\,i} e_i \quad a_{1\cdot s\cdot 1} = a_{1\cdot s\cdot 1}^{\,i} e_i \\ a_2^{ij} = a_{2\cdot s\cdot 0}^{\,i}\, a_{2\cdot s\cdot 1}^{\,j} \quad a_{2\cdot s\cdot 0} = a_{2\cdot s\cdot 0}^{\,i} e_i \quad a_{2\cdot s\cdot 1} = a_{2\cdot s\cdot 1}^{\,i} e_i \end{cases}$$

**Theorem 5.1.** *Equation* (5.2) *is equivalent to equation*

(5.3) $$x^k a_k^r = b^r$$

(5.4) $$a_k^r = a_1^{ij} C_{kj}^p C_{ip}^r + a_2^{ij} C_{ik}^p C_{pj}^r$$

*Proof.* From the equation (3.1), it follows that

(5.5) $$\overline{e}_i(x\overline{e}_j) = x^k \overline{e}_i(\overline{e}_k \overline{e}_j) = x^k \overline{e}_i C_{kj}^p \overline{e}_p = x^k C_{kj}^p C_{ip}^r \overline{e}_r$$

(5.6) $$(\overline{e}_i x)\overline{e}_j = x^k (\overline{e}_i \overline{e}_k)\overline{e}_j = x^k C_{ik}^p \overline{e}_p \overline{e}_j = x^k C_{ik}^p C_{pj}^r \overline{e}_r$$

Equations (5.3), (5.4) follow from equations (5.5), (5.6). □

**Theorem 5.2.** *If determinant*

(5.7) $$\det \|a_i^j\| \neq 0$$

*then equation* (5.1) *has only one solution.*

*If determinant equal* 0, *then $F$-linear dependence of vector $b$ from vectors $a_i^j \overline{e}_j$ is condition of existence of solution. In this case, equation* (4.1) *has infinitely many solutions. Otherwise equation does not have solution.*

*Proof.* The statement of the theorem is corollary of the theory of linear equations over field. □

## 6. Invertibility of Elements in $F$-Algebra

**Theorem 6.1.** $a \in A$ *does not have right inverse iff*

(6.1) $$\det \|C_{ij}^k a^i\| = 0$$

*Proof.* If $a \in A$ has right inverse, then the equation

(6.2) $$ax = \overline{e}_0$$

has solution. Using coordinates relative to the basis $\overline{\overline{e}}$, we can write the equation (6.2) in following form

(6.3) $$C_{ij}^k a^i x^j = \delta_0^k$$

Since number of variables in the system of linear equations (6.3) is equal to the number of equations, then $a \in A$ has right inverse iff

(6.4) $$\det \|C_{ij}^k a^i\| \neq 0$$

Therefore, the condition (6.1) is equivalent to the irreversibility of $a$ on the right. $\square$

**Theorem 6.2.** $a \in A$ is left zero divisor iff

(6.5) $$\det \|C_{ij}^k a^i\| = 0$$

*Proof.* According to the definition [1]-10.17, page 342, $a \in A$ is left zero divisor, if $a \neq 0$ and there exists

(6.6) $$x \neq 0$$

such that

(6.7) $$ax = 0$$

Using coordinates relative to the basis $\overline{\overline{e}}$, we can write the equation (6.7) in following form

(6.8) $$C_{ij}^k a^i x^j = 0$$

Since number of variables in the system of linear equations (6.8) is equal to the number of equations, then the condition (6.5) is equivalent to the condition (6.6). $\square$

**Theorem 6.3.** $a \in A$ is left zero divisor iff $a \in A$ does not have right inverse.

*Proof.* The theorem follows from theorems 6.1 and 6.2. $\square$

**Example 6.4.** Basis of algebra of $n \times n$ matrices consists of matrices

$$\overline{e}_k^i = (\delta_j^i \delta_k^l)$$

We define the product of matrices $\overline{e}_q^p$ and $\overline{e}_t^s$ by equation

(6.9) $$\overline{e}_q^p \overline{e}_t^s = (\delta_j^p \delta_q^l \delta_m^s \delta_t^j) = (\delta_t^p \delta_q^l \delta_m^s) = \delta_t^p \overline{e}_q^s = \delta_t^p \delta_i^s \delta_q^j \overline{e}_j^i$$

From the equation (6.9) it follows that

(6.10) $$C_{i \cdot q \cdot t}^{\cdot j \; p \; s} = \delta_t^p \delta_i^s \delta_q^j$$

From the equation (6.10) and theorem 6.1 it follows that the matrix $a$ does not have right inverse, if

(6.11) $$\det \|C_{i \cdot q \cdot t}^{\cdot j \; p \; s} a_s^t\| = \det \|\delta_t^p \delta_i^s \delta_q^j a_s^t\| = \det \|\delta_q^j a_i^p\| = 0$$

If we enumerate rows and columns of matrix (6.11) by index $^i_j$, then we will see that matrix (6.11) is matrix

$$(6.12) \qquad < A, E_n >= \begin{pmatrix} a^1_1 E_n & ... & a^1_n E_n \\ ... & ... & ... \\ a^n_1 E_n & ... & a^n_n 1 E_n \end{pmatrix}$$

where $E_n$ is $n \times n$ identity matrix.

Let us prove by induction over $m$ that

$$(6.13) \qquad \det < A, E_m >= (\det A)^m$$

For $m = 1$, $< A, E_1 >= A$. Therefore,

$$\det < A, E_1 >= \det A$$

Let the statement is true for $m - 1$. We can represent minor $a^i_j E_m$ in the form

$$\begin{pmatrix} a^i_j & 0 \\ 0 & a^i_j E_{m-1} \end{pmatrix}$$

We select in each minor $a^i_j E_m$ first element in first row, we get minor

$$\begin{pmatrix} a^1_1 & ... & a^1_n \\ ... & ... & ... \\ a^n_1 & ... & a^n_n \end{pmatrix} = A$$

that has algebraic complement

$$\begin{pmatrix} a^1_1 E_{m-1} & ... & a^1_n E_{m-1} \\ ... & ... & ... \\ a^n_1 E_{m-1} & ... & a^n_n E_{m-1} \end{pmatrix} =< A, E_{m-1} >$$

The other minors in the corresponding set of columns have at least one row consisting from 0. According to Laplace expansion theorem ([2], p. 51, [6], p. 259)

$$\det < A, E_m >= \det A \det < A, E_{m-1} >= \det A (\det A)^{m-1} = (\det A)^m$$

Therefore, from theorem 6.1 it follows that matrix $A$ does not have right inverse iff $(\det A)^m = 0$. It is evident that this condition is equivalent to request $\det A = 0$. $\square$

## 7. EQUATION $ax - xa = 1$

**Theorem 7.1.** *Equation*

$$ax - xa = 1$$

*in quaternion algebra does not have solutions.*

*Proof.* From equation, it follows that

$$(ax)_0 - (xa)_0 = 1$$

However, in quaternion algebra, it is true that

$$(ax)_0 = (xa)_0$$

□

**Theorem 7.2.** *Equation*
$$ax - xa = 1$$
*in algebra of matrix of order $2$ does not have solutions.*

*Proof.* The theorem immediately follows from calculations

$$\begin{pmatrix} a_1^1 & a_2^1 \\ a_1^2 & a_2^2 \end{pmatrix} \begin{pmatrix} x_1^1 & x_2^1 \\ x_1^2 & x_2^2 \end{pmatrix} - \begin{pmatrix} x_1^1 & x_2^1 \\ x_1^2 & x_2^2 \end{pmatrix} \begin{pmatrix} a_1^1 & a_2^1 \\ a_1^2 & a_2^2 \end{pmatrix} = \begin{pmatrix} 1 & 0 \\ 0 & 1 \end{pmatrix}$$

$$\begin{pmatrix} a_1^1 x_1^1 + a_2^1 x_1^2 & a_1^1 x_2^1 + a_2^1 x_2^2 \\ a_1^2 x_1^1 + a_2^2 x_1^2 & a_1^2 x_2^1 + a_2^2 x_2^2 \end{pmatrix} - \begin{pmatrix} x_1^1 a_1^1 + x_2^1 a_1^2 & x_1^1 a_2^1 + x_2^1 a_2^2 \\ x_1^2 a_1^1 + x_2^2 a_1^2 & x_1^2 a_2^1 + x_2^2 a_2^2 \end{pmatrix} = \begin{pmatrix} 1 & 0 \\ 0 & 1 \end{pmatrix}$$

$$\begin{pmatrix} a_2^1 x_1^2 - a_1^2 x_2^1 & a_1^1 x_2^1 + a_2^1 x_2^2 - a_2^1 x_1^1 - a_2^2 x_2^1 \\ a_1^2 x_1^1 + a_2^2 x_1^2 - a_1^1 x_1^2 - a_1^2 x_2^2 & a_1^2 x_2^1 - a_2^1 x_1^2 \end{pmatrix} = \begin{pmatrix} 1 & 0 \\ 0 & 1 \end{pmatrix}$$

$$\begin{cases} -a_1^2 x_2^1 & +a_2^1 x_1^2 & = 1 \\ -a_2^1 x_1^1 & +(a_1^1 - a_2^2) x_2^1 & +a_2^1 x_2^2 & = 0 \\ +a_1^2 x_1^1 & +(a_2^2 - a_1^1) x_1^2 & -a_1^2 x_2^2 & = 0 \\ +a_1^2 x_2^1 & -a_2^1 x_1^2 & = 1 \end{cases}$$

It is evident that first and fourth equation are incompatible.                □

**Theorem 7.3.** *The equation*
$$ax - xa = 1$$
*in the algebra A has only one solution if*
$$\det \|(C_{ij}^k - C_{ji}^k) a^i\| \neq 0$$

*Proof.* The theorem is corollary of the theorem 4.1.                □

# Линейное уравнение в конечномерной алгебре

Александр Клейн

Аннотация. В статье рассмотрены методы решения уравнения вида

$$a_{s\cdot 0}\, x\, a_{s\cdot 1} = b$$

в алгебре, конечномерной над полем.

Содержание

## 1. Предисловие

Мои недавние исследования в области математического анализа ([4]) и аффинной геометрии над телом ([5]) привели меня к необходимости решать линейные уравнения вида

$$(1.1) \qquad a_{s\cdot 0}\, x\, a_{s\cdot 1} = b$$

или системы подобных уравнений over division ring. Основная задача, которая меня интересует, - это найти обратное преобразование, так как это основная операция при при преобразовании тензора, а так же при поднятии и опускании индекса в эвклидовом пространстве.

В этой статье я исследовал возможность решения простейшего уравнения. В современной математике и физике рассматриваются структуры где операция произведения может быть неассоциативной. Поэтому меня интересует, какие результаты сохранятся, если я буду рассматривать линейную алгебру в конечномерной алгебре $A$ над полем.

Коэффициенты уравнения (1.1) принадлежат тензорному произведению $A \otimes A$. Это позволяет применить разработанные ранее методы к решению уравнения (1.1). В статье я рассматриваю два метода решения уравнения (1.1).

Aleks_Kleyn@MailAPS.org.
http://sites.google.com/site/AleksKleyn/.
http://arxiv.org/a/kleyn_a_1.
http://AleksKleyn.blogspot.com/.

Сперва я рассматриваю алгебру $A$ как векторное пространство над полем $F$. Это позволяет заменить уравнение (1.1) системой линейных уравнений, которые мы умеем решать. Сразу становится очевидным, что уравнение (1.1) может иметь один корень, бесконечно много корней либо ни одного.

Хотя стандартное представление линейного выражения не всегда самая компактная форма записи, оно позволяет упростить построение. Второй метод решения уравнения (1.1) опирается на стандартное представление линейного выражения. Когда уравнение имеет единственное решение, этот метод позволяет найти обратное отображение (задача, которая важна для меня в тензорном исчислении), а так же лучше понять, что значит вырожденное линейное отображение (это необходимо для построения гомологий алгебры; в частности, что бы понять как может выглядеть уравнение Коши-Римана в алгебре и какие функции ему удовлетворяют).

Выражение (4.8) как запись тождественного преобразования выглядит необычно. В первый раз я с таким выражением столкнулся в равенстве [3]-(6.4.3). Но это было самое начало моего исследования в этой области, и я не был готов в то время понять всей глубины мною увиденного.

## 2. Соглашения

**Соглашение 2.1.** *В любом выражении, где появляется индекс, я предполагаю, что этот индекс может иметь внутреннюю структуру. Например, при рассмотрении алгебры $A$ координаты $a \in A$ относительно базиса $\overline{\overline{e}}$ пронумерованы индексом $i$. Это означает, что $a$ является вектором. Однако, если $a$ является матрицей, нам необходимо два индекса, один нумерует строки, другой - столбцы. В том случае, когда мы уточняем структуру индекса, мы будем начинать индекс с символа $\cdot$ в соответствующей позиции. Например, если я рассматриваю матрицу $a^i_j$ как элемент векторного пространства, то я могу записать элемент матрицы в виде $a^{\cdot i}_{\cdot j}$.* □

**Соглашение 2.2.** *В выражении вида*

$$a_{s \cdot 0} x a_{s \cdot 1}$$

*предполагается сумма по индексу $s$.* □

**Соглашение 2.3.** *Для данного поля $F$, если не оговорено противное, мы будем рассматривать конечномерную $F$-алгебру.* □

**Соглашение 2.4.** *Пусть $A$ - свободная конечно мерная алгебра. При разложении элемента алгебры $A$ относительно базиса $\overline{\overline{e}}$ мы пользуемся одной и той же корневой буквой для обозначения этого элемента и его координат. Однако в алгебре не принято использовать векторные обозначения. В выражении $a^2$ не ясно - это компонента разложения элемента $a$ относительно базиса или это операция возведения в степень. Для облегчения чтения текста мы будем индекс элемента алгебры выделять цветом. Например,*

$$a = a^{i} \overline{e}_{i}$$

□

**Соглашение 2.5.** *Если свободная конечномерная алгебра имеет единицу, то мы будем отождествлять вектор базиса $\overline{e}_0$ с единицей алгебры.* □

Без сомнения, у читателя могут быть вопросы, замечания, возражения. Я буду признателен любому отзыву.

## 3. Произведение в алгебре

Пусть $F$ - поле и $A$ - конечномерная алгебра над полем $F$. Пусть $\overline{\overline{e}}$ - базис алгебры $A$ над полем $F$. Произведение базисных векторов определено согласно правилу

$$\overline{e}_i \overline{e}_j = C_{ij}^k \overline{e}_k \tag{3.1}$$

где $C_{ij}^k$ - структурные константы алгебры $A$ над полем $F$. Поскольку операция не предполагается ни коммутативной, ни ассоциативной, мы не накладываем никаких ограничений на структурные константы.

Из равенства (3.1) следует, что произведение $a, b \in A$ можно получить согласно правилу

$$(ab)^k = C_{ij}^k a^i b^j \tag{3.2}$$

## 4. Линейное уравнение в ассоциативной алгебре

Линейное уравнение в ассоциативной алгебре $A$ мы будем записывать в виде

$$a_{s\cdot 0}\, x\, a_{s\cdot 1} = b \tag{4.1}$$

Согласно теореме [3]-9.1.9 мы можем записать уравнение (4.1) в стандартной форме

$$a^{ij}\overline{e}_i x \overline{e}_j = b$$
$$a^{ij} = a_{s\cdot 0}^{\ i} a_{s\cdot 1}^{\ j} \quad a_{s\cdot 0} = a_{s\cdot 0}^{\ i} e_i \quad a_{s\cdot 1} = a_{s\cdot 1}^{\ i} e_i \tag{4.2}$$

Согласно теореме [3]-9.1.10 уравнение (4.2) эквивалентно уравнению

$$x^i a_i^j = b^j \tag{4.3}$$

$$a_i^j = a^{kr} C_{ki}^p C_{pr}^j \tag{4.4}$$

**Теорема 4.1.** *Если определитель*

$$\det \|a_i^j\| \neq 0 \tag{4.5}$$

*то уравнение* (4.1) *имеет единственное решение.*

*Если определитель равен* $0$, *то условием существования решения является $F$-линейная зависимость вектора $b$ от векторов $a_i^j \overline{e}_j$. В этом случае, уравнение* (4.1) *имеет бесконечно много решений. В противном случае уравнение не имеет решений.*

*Доказательство.* Утверждение теоремы является следствием теории линейных уравнений над полем. $\square$

**Теорема 4.2.** *Пусть уравнение* (4.2) *удовлетворяет условию* (4.5). *Если равенство* (4.2) *рассматривать как преобразование алгебры $A$, то обратное преобразование можно записать в виде*

$$x = c^{pq} \overline{e}_p b \overline{e}_q \tag{4.6}$$

*где компоненты $c^{pq}$ удовлетворяют уравнению*

$$\delta_0^r \delta_0^s = a^{ij} c^{pq} C_{ip}^r C_{qj}^s \tag{4.7}$$

*Доказательство.* Согласно теореме 4.1 решение уравнения (4.2) единственно. Так как $x$ является линейной функцией от $b$, то мы рассматриваем стандартное представление (4.6) этой функции. Из равенств (4.2), (4.6) следует

$$(4.8) \qquad b = a^{ij} c^{pq} \overline{e}_i \overline{e}_p b \overline{e}_q \overline{e}_j = a^{ij} c^{pq} C_{ip}^r C_{qj}^s \overline{e}_r b \overline{e}_s$$

Поскольку $b$ - произвольно, равенство (4.8) является записью тождественного преобразования. Равенство (4.7) является следствием равенства (4.8) и теоремы [3]-9.1.15. □

## 5. Линейное уравнение в неассоциативной алгебре

Линейное уравнение в неассоциативной алгебре $A$ мы будем записывать в виде

$$(5.1) \qquad a_{1 \cdot s \cdot 0} (x \, a_{1 \cdot s \cdot 1}) + (a_{2 \cdot t \cdot 0} \, x) \, a_{2 \cdot t \cdot 1} = b$$

Согласно теореме [3]-9.1.9 мы можем записать уравнение (5.1) в стандартной форме

$$(5.2) \qquad \begin{cases} a_1^{ij} \overline{e}_i (x \overline{e}_j) + a_2^{ij} (\overline{e}_i x) \overline{e}_j = b \\ a_1^{ij} = a_{1 \cdot s \cdot 0}^{\ i} a_{1 \cdot s \cdot 1}^{\ j} \quad a_{1 \cdot s \cdot 0} = a_{1 \cdot s \cdot 0}^{\ i} e_i \quad a_{1 \cdot s \cdot 1} = a_{1 \cdot s \cdot 1}^{\ i} e_i \\ a_2^{ij} = a_{2 \cdot s \cdot 0}^{\ i} a_{2 \cdot s \cdot 1}^{\ j} \quad a_{2 \cdot s \cdot 0} = a_{2 \cdot s \cdot 0}^{\ i} e_i \quad a_{2 \cdot s \cdot 1} = a_{2 \cdot s \cdot 1}^{\ i} e_i \end{cases}$$

**Теорема 5.1.** *Уравнение* (5.2) *эквивалентно уравнению*

$$(5.3) \qquad x^k a_k^r = b^r$$

$$(5.4) \qquad a_k^r = a_1^{ij} C_{kj}^p C_{ip}^r + a_2^{ij} C_{ik}^p C_{pj}^r$$

*Доказательство.* Из равенства (3.1) следует

$$(5.5) \qquad \overline{e}_i (x \overline{e}_j) = x^k \overline{e}_i (\overline{e}_k \overline{e}_j) = x^k \overline{e}_i C_{kj}^p \overline{e}_p = x^k C_{kj}^p C_{ip}^r \overline{e}_r$$

$$(5.6) \qquad (\overline{e}_i x) \overline{e}_j = x^k (\overline{e}_i \overline{e}_k) \overline{e}_j = x^k C_{ik}^p \overline{e}_p \overline{e}_j = x^k C_{ik}^p C_{pj}^r \overline{e}_r$$

Равенства (5.3), (5.4) следуют из равенств (5.5), (5.6). □

**Теорема 5.2.** *Если определитель*

$$(5.7) \qquad \det \|a_i^j\| \neq 0$$

*то уравнение* (5.1) *имеет единственное решение.*

*Если определитель равен $0$, то условием существования решения является $F$-линейная зависимость вектора $b$ от векторов $a_i^j \overline{e}_j$ . В этом случае, уравнение* (4.1) *имеет бесконечно много решений. В противном случае уравнение не имеет решений.*

*Доказательство.* Утверждение теоремы является следствием теории линейных уравнений над полем. □

## 6. Обратимость элемента в $F$-алгебре

**Теорема 6.1.** *$a \in A$ не имеет правого обратного тогда и только тогда, когда*

$$\det \|C_{ij}^k a^i\| = 0 \tag{6.1}$$

*Доказательство.* Если $a \in A$ имеет правый обратный, то уравнение

$$ax = \overline{e}_0 \tag{6.2}$$

имеет решение. В координатах относительно базиса $\overline{\overline{e}}$, уравнение (6.2) имеет вид

$$C_{ij}^k a^i x^j = \delta_0^k \tag{6.3}$$

Так как число неизвестных в системе линейных уравнений (6.3) равно числу уравнений, то $a \in A$ имеет правый обратный тогда и только тогда, когда

$$\det \|C_{ij}^k a^i\| \neq 0 \tag{6.4}$$

Следовательно, условие (6.1) эквивалентно необратимости $a$ справа. □

**Теорема 6.2.** *$a \in A$ является левым делителем нуля тогда и только тогда, когда*

$$\det \|C_{ij}^k a^i\| = 0 \tag{6.5}$$

*Доказательство.* Согласно определению [1]-10.17, страница 342, $a \in A$ является левым делителем нуля, если $a \neq 0$ и существует

$$x \neq 0 \tag{6.6}$$

такое, что

$$ax = 0 \tag{6.7}$$

В координатах относительно базиса $\overline{\overline{e}}$, уравнение (6.7) имеет вид

$$C_{ij}^k a^i x^j = 0 \tag{6.8}$$

Так как число неизвестных в системе линейных уравнений (6.8) равно числу уравнений, то условие (6.5) эквивалентно условию (6.6). □

**Теорема 6.3.** *$a \in A$ является левым делителем нуля тогда и только тогда, когда $a \in A$ не имеет правого обратного.*

*Доказательство.* Утверждение теоремы является следствием теорем 6.1 и 6.2. □

**Пример 6.4.** Базис в алгебре $n \times n$ матриц состоит из матриц вида

$$\overline{e}_k^i = (\delta_j^i \delta_k^l)$$

Произведение матриц $\overline{e}_q^p$ и $\overline{e}_t^s$ определенно равенством

$$\overline{e}_q^p \overline{e}_t^s = (\delta_j^p \delta_q^l \delta_m^s \delta_t^j) = (\delta_t^p \delta_q^l \delta_m^s) = \delta_t^p \overline{e}_q^s = \delta_t^p \delta_i^s \delta_q^j \overline{e}_j^i \tag{6.9}$$

Из равенства (6.9) следует

$$C^{\cdot j\ p\ s}_{i\ \cdot q\ \cdot t} = \delta_t^p \delta_i^s \delta_q^j \tag{6.10}$$

Из равенства (6.10) и теоремы 6.1 следует, что матрица $a$ необратима справа, если

$$\det \|C^{\cdot j\ p\ s}_{i\ \cdot q\ \cdot t} a_s^t\| = \det \|\delta_t^p \delta_i^s \delta_q^j a_s^t\| = \det \|\delta_q^j a_i^p\| = 0 \tag{6.11}$$

Если мы пронумеруем строки и столбцы матрицы (6.11) индексом $^i_j$, то мы увидим, что матрица (6.11) является матрицей вида

(6.12) $$<A, E_n> = \begin{pmatrix} a^1_1 E_n & ... & a^1_n E_n \\ ... & ... & ... \\ a^n_1 E_n & ... & a^n_n 1 E_n \end{pmatrix}$$

где $E_n$ - $n \times n$ единичная матрица.

Докажем индукцией по $m$, что

(6.13) $$\det <A, E_m> = (\det A)^m$$

Для $m = 1$, $<A, E_1> = A$. Следовательно,

$$\det <A, E_1> = \det A$$

Пусть утверждение верно для $m-1$. Минор $a^i_j E_m$ можно представить в виде

$$\begin{pmatrix} a^i_j & 0 \\ 0 & a^i_j E_{m-1} \end{pmatrix}$$

Выбрав в каждом миноре $a^i_j E_m$ первый элемент в первой сторке, мы получим минор

$$\begin{pmatrix} a^1_1 & ... & a^1_n \\ ... & ... & ... \\ a^n_1 & ... & a^n_n \end{pmatrix} = A$$

имеющий алгебраическое дополнение

$$\begin{pmatrix} a^1_1 E_{m-1} & ... & a^1_n E_{m-1} \\ ... & ... & ... \\ a^n_1 E_{m-1} & ... & a^n_n E_{m-1} \end{pmatrix} = <A, E_{m-1}>$$

Остальные миноры в соответствующем множестве столбцов имеют, по крайней мере, одну строку, состоящую из 0. Согласно теореме Лапласа ([2], с. 51, [6], с. 259)

$$\det <A, E_m> = \det A \det <A, E_{m-1}> = \det A (\det A)^{m-1} = (\det A)^m$$

Следовательно, из теоремы 6.1 следует, что матрица $A$ необратима справа тогда и только тогда, когда $(\det A)^m = 0$. Очевидно, что это условие эквивалентно требованию $\det A = 0$. □

## 7. Уравнение $ax - xa = 1$

**Теорема 7.1.** *Уравнение*

$$ax - xa = 1$$

*в алгебре кватернионов не имеет решений.*

*Доказательство.* Из уравнения следует
$$(ax)_{\mathbf{0}} - (xa)_{\mathbf{0}} = 1$$
Но в алгебре кватернионов
$$(ax)_{\mathbf{0}} = (xa)_{\mathbf{0}}$$
□

**Теорема 7.2.** *Уравнение*
$$ax - xa = 1$$
*в алгебре матриц порядка 2 не имеет решений.*

*Доказательство.* Теорема непосредственно следует из вычислений
$$\begin{pmatrix} a_1^1 & a_2^1 \\ a_1^2 & a_2^2 \end{pmatrix} \begin{pmatrix} x_1^1 & x_2^1 \\ x_1^2 & x_2^2 \end{pmatrix} - \begin{pmatrix} x_1^1 & x_2^1 \\ x_1^2 & x_2^2 \end{pmatrix} \begin{pmatrix} a_1^1 & a_2^1 \\ a_1^2 & a_2^2 \end{pmatrix} = \begin{pmatrix} 1 & 0 \\ 0 & 1 \end{pmatrix}$$

$$\begin{pmatrix} a_1^1 x_1^1 + a_2^1 x_1^2 & a_1^1 x_2^1 + a_2^1 x_2^2 \\ a_1^2 x_1^1 + a_2^2 x_1^2 & a_1^2 x_2^1 + a_2^2 x_2^2 \end{pmatrix} - \begin{pmatrix} x_1^1 a_1^1 + x_2^1 a_1^2 & x_1^1 a_2^1 + x_2^1 a_2^2 \\ x_1^2 a_1^1 + x_2^2 a_1^2 & x_1^2 a_2^1 + x_2^2 a_2^2 \end{pmatrix} = \begin{pmatrix} 1 & 0 \\ 0 & 1 \end{pmatrix}$$

$$\begin{pmatrix} a_2^1 x_1^2 - a_1^2 x_2^1 & a_1^1 x_2^1 + a_2^1 x_2^2 - a_2^1 x_1^1 - a_2^2 x_2^1 \\ a_1^2 x_1^1 + a_2^2 x_1^2 - a_1^1 x_1^2 - a_1^2 x_2^2 & a_1^2 x_2^1 - a_2^1 x_1^2 \end{pmatrix} = \begin{pmatrix} 1 & 0 \\ 0 & 1 \end{pmatrix}$$

$$\begin{cases} -a_1^2 x_2^1 & +a_2^1 x_1^2 & = 1 \\ -a_2^1 x_1^1 & +(a_1^1 - a_2^2) x_2^1 & +a_2^1 x_2^2 & = 0 \\ +a_1^2 x_1^1 & +(a_2^2 - a_1^1) x_1^2 & -a_1^2 x_2^2 & = 0 \\ +a_1^2 x_2^1 & -a_2^1 x_1^2 & = 1 \end{cases}$$

Очевидно, первое и четвёртое уравнение несовместимы. □

**Теорема 7.3.** *Уравнение*
$$ax - xa = 1$$
*в алгебре A имеет единственное решение при условии*
$$\det \|(C_{ij}^{\boldsymbol{k}} - C_{ji}^{\boldsymbol{k}}) a^{\boldsymbol{i}}\| \neq 0$$

*Доказательство.* Теорема является следствием теоремы 4.1. □

## 8. Список литературы